\documentclass[12pt,a4paper,draft]{amsart}

\usepackage{amssymb,amsrefs}
\usepackage[all]{xy}

\allowdisplaybreaks

\renewcommand{\bf}{\mathbf{f}}

\newcommand{\bd}{\mathbf{d}}
\newcommand{\be}{\mathbf{e}}
\newcommand{\bu}{\mathbf{u}}
\newcommand{\bv}{\mathbf{v}}
\newcommand{\bw}{\mathbf{w}}
\newcommand{\bx}{\mathbf{x}}
\newcommand{\by}{\mathbf{y}}

\newcommand{\bbA}{\mathbb{A}}
\newcommand{\bbD}{\mathbb{D}}
\newcommand{\bbM}{\mathbb{M}}
\newcommand{\bbN}{\mathbb{N}}
\newcommand{\bbZ}{\mathbb{Z}}

\newcommand{\calC}{\mathcal{C}}
\newcommand{\calF}{\mathcal{F}}
\newcommand{\calO}{\mathcal{O}}
\newcommand{\calU}{\mathcal{U}}
\newcommand{\calZ}{\mathcal{Z}}

\let\mod=\undefined
\DeclareMathOperator{\GL}{GL} %
\DeclareMathOperator{\tr}{tr} %
\DeclareMathOperator{\ind}{ind} %
\DeclareMathOperator{\ini}{in} %
\DeclareMathOperator{\Ker}{Ker} %
\DeclareMathOperator{\mod}{mod} %
\DeclareMathOperator{\rep}{rep} %
\DeclareMathOperator{\Stab}{Stab} %

\newcommand{\ol}{\overline}

\newtheorem{coro}{Corollary}[section]
\newtheorem{lemm}[coro]{Lemma}
\newtheorem{prop}[coro]{Proposition}
\newtheorem{theo}[coro]{Theorem}

\title[Orbit closures for directing modules]{Normality of
orbit closures for directing modules over tame algebras}

\author{Grzegorz Bobi\'nski}

\address{Faculty of Mathematics and Computer Science \\
Nicolaus Coperincus University \\ Chopina 12/18 \\ 87-100 Toru\'n
\\ Poland}

\email{gregbob@mat.uni.torun.pl}

\address{Faculty of Mathematics and Computer Science \\
Nicolaus Coperincus University \\ Chopina 12/18 \\ 87-100 Toru\'n
\\ Poland}

\email{gzwara@mat.uni.torun.pl}



\author{Grzegorz Zwara}

\date{\today}

\begin{document}

\begin{abstract}
We show that the orbit closures for directing modules over tame
algebras are normal and Cohen--Macaulay. The proof is based on
deformations to normal toric varieties.
\end{abstract}

\maketitle

\section{Introduction and the main results}

Throughout the paper $k$ denotes a fixed algebraically closed
field. By an algebra we mean an associative $k$-algebra with
identity, and by a module a finite dimensional left module.
Furthermore, for an algebra $A$, $\mod A$ stands for the category
of finite dimensional left $A$-modules. By $\bbN$ and $\bbZ$ we
denote the sets of nonnegative integers and integers,
respectively. Finally, if $i$ and $j$ are integers, then by $[i,
j]$ we denote the set of all integers $k$ such that $i \leq k \leq
j$.

Let $d$ be a positive integer and denote by $\bbM (d)$ the algebra
of $d \times d$-matrices with coefficients in $k$. For an algebra
$A$ the set $\mod_A (d)$ of the $A$-module structures on the
vector space $k^d$ has a natural structure of an affine variety.
Indeed, if $A \simeq k \langle X_1, \ldots, X_t \rangle / I$ for
$t > 0$ and a two-sided ideal $I$, then $\mod_A (d)$ can be
identified with the closed subset of $(\bbM (d))^t$ given by the
vanishing of the entries of all matrices $\rho (X_1, \ldots, X_t)$
for $\rho \in I$. Moreover, the general linear group $\GL (d)$
acts on $\mod_A (d)$ by conjugations and the $\GL(d)$-orbits in
$\mod_A (d)$ correspond bijectively to the isomorphism classes of
$d$-dimensional left $A$-modules. We shall denote by $\calO_M$ the
$\GL (d)$-orbit in $\mod_A (d)$ corresponding to (the isomorphism
class of) a $d$-dimensional module $M$ in $\mod A$. It is an
interesting task to study geometric properties of the Zariski
closure $\ol{\calO}_M$ of $\calO_M$.

The above problem can also be formulated in terms of
representations of finite quivers instead of modules over
algebras. Here, by a finite quiver $\Sigma$ we mean a finite set
$\Sigma_0$ of vertices and a finite set $\Sigma_1$ of arrows
together with two maps $s, t : \Sigma_1 \to \Sigma_0$, which
assign to an arrow its starting and terminating vertex,
respectively. Let $\bd = (d_x)_{x \in \Sigma_0} \in
\bbN^{\Sigma_0}$ be a dimension vector and let $\bbM (m, n)$
denote the space of $m \times n$-matrices with coefficients in
$k$. The affine space
\[
\rep_\Sigma (\bd)= \prod_{\alpha \in \Sigma_1} \bbM (d_{t \alpha},
d_{s \alpha})
\]
is called a variety of representations of $\Sigma$. The product
$\GL (\bd) = \prod_{x \in \Sigma_0} \GL (d_x)$ of general linear
groups acts on $\rep_\Sigma (\bd)$ by conjugations:
\[
g \cdot V = (g_{t \alpha} V_\alpha g_{s \alpha}^{-1})_{\alpha \in
\Sigma_1}
\]
for $g = (g_x)_{x \in \Sigma_0} \in \GL (\bd)$ and $V =
(V_\alpha)_{\alpha \in \Sigma_1} \in \rep_\Sigma (\bd)$. The orbit
of $V \in \rep_\Sigma (\bd)$ with respect to this action is
denoted by $\calO_V$, and its closure by $\ol{\calO}_V$. In fact,
the module varieties and varieties of representations of quivers
are closely related to each other (see \cite{Bon1} for details).
In particular, for any algebra $A$ there is a uniquely determined
quiver $\Sigma$ (called the Gabriel quiver of $A$) such that for
each $d \geq 1$ and $M \in \mod_A (d)$ there are a dimension
vector $\bd \in \bbN^{\Sigma_0}$ and $V \in \rep_\Sigma (\bd)$
such that $\ol{\calO}_M$ is isomorphic to the associated fibre
bundle $\GL (d) \times_{\GL (\bd)} \ol{\calO}_V$. Hence
$\ol{\calO}_M$ is normal, Cohen-Macaulay, unibranch or regular in
some codimension if and only if $\ol{\calO}_V$ is.

The orbit closures are normal and Cohen--Macaulay varieties (with
rational singularities in characteristic zero) provided $\Sigma$
is a Dynkin quiver of type $\bbA_n$ or $\bbD_n$ (\cites{BobZw1,
BobZw2}), or $A$ is a Brauer tree algebra (\cite{SkZw}). Moreover,
they are regular in codimension one if $\Sigma$ is the Kronecker
quiver (\cite{BeBon}), or $A$ is a representation finite algebra
(\cite{Zw2}), i.e., a set $\ind A$ of chosen representatives of
isomorphism classes of indecomposable $A$-modules is finite.
Another result states that the variety $\ol{\calO}_M$ is unibranch
if there are only finitely many modules $U$ in $\ind A$ such that
there is a monomorphism from $U$ to $M^i$ for some $i > 0$
(\cite{Zw3}). On the other hand, there exists an orbit closure in
$\rep_\Sigma ((3,3))$, where $\Sigma$ is the Kronecker quiver,
which is neither unibranch nor Cohen--Macaulay (see \cite{Zw1}).

We say that an algebra $A$ is tame if we can chose $\ind A$ in
such a way that for every $d > 0$ all $d$-dimensional modules in
$\ind A$ can be described by finitely many one-parameter families.
According to Drozd's Tame and Wild Theorem (\cite{Dr}, see also
\cite{CB}) there is a chance to classify modules only for tame
algebras. An indecomposable module $M$ in $\mod A$ is called
directing if there exists no sequence
\[
M = M_0 \xrightarrow{f_1} M_1 \to \cdots \to M_{m - 1}
\xrightarrow{f_m} M_m = M
\]
in $\mod A$, where $m > 0$, $M_1$, \ldots, $M_{m - 1}$ belong to
$\ind A$ and $f_1$, \ldots, $f_m$ are nonzero nonisomorphisms.
Bongartz investigated from the geometric point of view a special
class of directing modules, so called preprojective ones
(see~\cite{Bon2}*{Proposition~6}). Further results in this
direction were obtained by Skowro\'nski and the first author
in~\cite{BobSk1} (see also~\cite{Bob} for the case of decomposable
directing modules). The main result of the paper is as follows.

\begin{theo} \label{theomain}
Let $M$ be an indecomposable directing module over a tame algebra.
Then the variety $\ol{\calO}_M$ is normal and Cohen-Macaulay.
\end{theo}

Using \cite{BobSk1}*{Theorem~2} (see
\cite{BobSk3}*{Proposition~2.4} for the correct list of algebras)
and the geometric equivalence described in~\cite{Bon1} we get that
$\ol{\calO}_M$ is isomorphic to the associated fibre bundle $\GL
(d) \times_{\GL (\bd)} \ol{\calO}_P$, where either $\ol{\calO}_P$
is a normal complete intersection, or up to duality, $P$ is
defined as follows. Let $0 \leq p \leq q \leq r \leq s \leq t$,
let $\Delta$ be the quiver
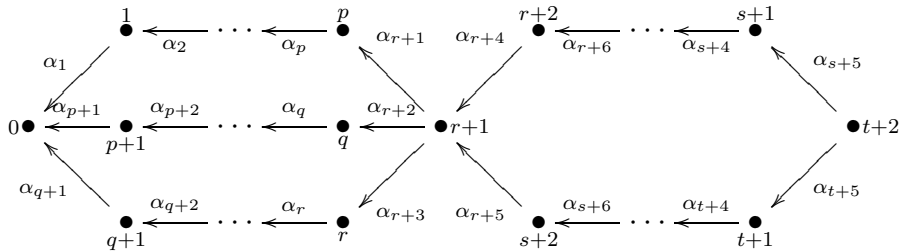
\[
\xymatrix{& \bullet \save*+!D{\scriptstyle 1} \restore
\ar[ld]_{\alpha_1} & \cdots \ar[l]^-{\alpha_2} & \bullet
\save*+!D{\scriptstyle p} \restore \ar[l]^-{\alpha_p} & & \bullet
\save*+!D{\scriptstyle r + 2} \restore \ar[ld]_(0.25){\alpha_{r +
4}} & \cdots \ar[l]^-{\alpha_{r + 6}} & \bullet
\save*+!D{\scriptstyle s + 1} \restore \ar[l]^-{\alpha_{s + 4}} \\
\bullet \save*+!R{\scriptstyle 0} \restore & \bullet
\save*+!U{\scriptstyle p + 1} \restore \ar[l]_{\alpha_{p + 1}} &
\cdots \ar[l]_-{\alpha_{p + 2}} & \bullet \save*+!U{\scriptstyle
q} \restore \ar[l]_-{\alpha_q} & \bullet \save*+!L{\scriptstyle r
+ 1} \restore \ar[lu]_(0.75){\alpha_{r + 1}} \ar[l]_{\alpha_{r +
2}} \ar[ld]^(0.75){\alpha_{r + 3}} & & & & \bullet
\save*+!L{\scriptstyle t + 2} \restore \ar[lu]_{\alpha_{s + 5}}
\ar[ld]^{\alpha_{t + 5}} \\ & \bullet \save*+!U{\scriptstyle q +
1} \restore \ar[lu]^{\alpha_{q + 1}} & \cdots \ar[l]_-{\alpha_{q +
2}} & \bullet \save*+!U{\scriptstyle r} \restore
\ar[l]_-{\alpha_r} & & \bullet \save*+!U{\scriptstyle s + 2}
\restore \ar[lu]^(0.25){\alpha_{r + 5}} & \cdots
\ar[l]_-{\alpha_{s + 6}} & \bullet \save*+!U{\scriptstyle t + 1}
\restore \ar[l]_-{\alpha_{t + 4}}}
\]
(if some of the inequalities between $0$, $p$, $q$, $r$, $s$ and
$t$ are equalities, then we obtain the obvious degenerated version
of the above quiver; see also a more detailed discussion about the
definition of the quiver $Q (p, q, r, s, t)$ after
Proposition~\ref{propideal} in Section~\ref{sect2}) and let $\bd$
be the dimension vector in $\bbN^{\Delta_0}$, whose $(r + 1)$th
coordinate equals $2$ and the remaining coordinates are $1$. Then
$P = P (p, q, r, s, t)$ is the point $(P_\alpha)_{\alpha \in
\Delta_1} \in \rep_\Delta (\bd)$ such that
\begin{gather*}
P_{\alpha_{r + 1}} = [
\begin{matrix}
1 & 0
\end{matrix}
], \quad P_{\alpha_{r + 2}} = [
\begin{matrix}
-1 & -1
\end{matrix}
], \quad P_{\alpha_{r + 3}} = [
\begin{matrix}
0 & 1
\end{matrix}
],
\\ %
P_{\alpha_{r + 4}} = [
\begin{matrix}
0 & 1
\end{matrix}
]^{\tr}, \quad P_{\alpha_{r + 5}} = [
\begin{matrix}
1 & 0
\end{matrix}
]^{\tr},
\end{gather*}
and the remaining matrices $P_\alpha$ are equal to $[1]$. Hence
Theorem~\ref{theomain} is a consequence of the following result.

\begin{theo} \label{theospec}
Let $P = P (p, q, r, s, t)$ for some integers $0 \leq p \leq q
\leq r \leq s \leq t$. Then the variety $\ol{\calO}_P$ is normal,
Cohen--Macaulay, and has rational singularities in characteristic
zero.
\end{theo}

The idea of the proof is to deform such varieties to toric normal
varieties using the so-called Sagbi-bases (see \cites{RoSw,
CoHeVa}). These normal toric varieties appear in the following
theorem.

\begin{theo} \label{theotoric}
Let $Q$ be a finite quiver without oriented cycles, let $\bd$ be
the dimension vector in $\bbN^{Q_0}$ with the coordinates equal to
$1$ and let $V$ be the point of $\rep_Q (\bd)$ given by the
matrices equal to $[1]$. Then $\ol{\calO}_V$ is a normal
toric variety.
\end{theo}

The paper is organized as follows. In Section~\ref{sect2} we prove
Theorem~\ref{theotoric} and investigate the equations defining the
toric varieties described in the theorem. Section~\ref{sect3} is
devoted to the proof of Theorem~\ref{theospec}.

The paper was written during the authors' joint stay at the
University of Berne. Authors gratefully acknowledge the support
from the Schweizerischer Nationalfonds and the Polish Scientific
Grant KBN No.~1 P03A 018 27.

\section{Toric varieties} \label{sect2}

Let $Q$ be a finite quiver without oriented cycles and let $\bd =
(d_i)_{i \in Q_0}$ be the dimension vector in $\bbN^{Q_0}$ with
all $d_i$ equal to $1$. Then the algebraic group $\GL (\bd) =
\prod_{i \in Q_0} k^*$ is a torus and the orbit closures in
$\rep_Q(\bd)$ are affine toric varieties (here we do not assume
that toric varieties are normal). In particular, this holds for
the orbit closure $\ol{\calO}_V$, where $V = (V_\alpha)_{\alpha
\in Q_1}$ is the point of $\rep_Q (\bd)$ with $V_\alpha = [1]$ for
any arrow $\alpha \in Q_1$. Let $\be_\alpha = \be_{t \alpha} -
\be_{s \alpha}$ for $\alpha \in Q_1$, where $(\be_i)_{i \in Q_0}$
is the standard basis of $\bbZ^{Q_0}$. It follows from the
definition of the action of $\GL(\bd)$ on $\rep_Q (\bd)$ that
$\ol{\calO}_V$ corresponds to the cone
\[
\calC_Q = \sum_{\alpha \in Q_1} \bbN \cdot \be_\alpha \subset
\bbZ^{Q_0},
\]
which means that the algebra $k [\ol{\calO}_V]$ of regular
functions on $\ol{\calO}_V$ may be identified with the subalgebra
of $k [T_i, T_i^{-1}]_{i \in Q_0}$ generated by $T^{\be_\alpha}$,
$\alpha\in Q_1$, where for $\bx = (x_i)_{i \in Q_0} \in
\bbZ^{Q_0}$ we put $T^{\bx} = \prod_{i \in Q_0} T_i^{x_i}$.
According to this identification, $k [\ol{\calO}_V]$ as a vector
space has a basis formed by $T^\bx$, $\bx \in \calC_Q$. It is
well-known that an affine toric variety is normal if and only if
the corresponding cone $\calC$ is saturated, i.e., a lattice point
$\bx$ belongs to $\calC$ whenever $\lambda \bx \in\calC$ for some\
$\lambda \in \bbN \setminus \{ 0 \}$.

For a vector $ \bx = (x_i)_{i \in Q_0} \in\bbZ^{Q_0}$ and a subset
$F$ of $Q_0$ we abbreviate by $\bx_F$ the sum $\sum_{i\in F} x_i$.
A subset $F$ of $Q_0$ is called a filter in $Q$ if
\[
s \alpha \in F \implies t \alpha \in F
\]
for any arrow $\alpha \in Q_1$. Let $X_Q$ be the subset of all
$\bx \in \bbZ^{Q_0}$ such that $\bx_{Q_0} = 0$ and $\bx_F \geq 0$
for any filter $F$ in $Q$. Obviously $X_Q$ is a saturated cone.
Hence Theorem~\ref{theotoric} is a consequence of the following
fact.

\begin{prop}
$\calC_Q = X_Q$.
\end{prop}

\begin{proof}
Obviously $\calC_Q \subseteq X_Q$. Let $\bx = (x_i)_{i \in Q_0}
\in X_Q$. In order to prove that $\bx \in \calC_Q$ we proceed by a
double induction, first: on the cardinality of $Q_0$, and second:
on the integer $\sum_{F \in \calF} \bx_F \geq 0$, where $\calF$ is
the set of all filters in $Q$.

Assume first that there is no arrow in $Q_1$ (for example, this
holds if $Q_0$ has only one element). Then for any $i \in Q_0$,
$\{ i \}$ is a filter in $Q$ and thus $x_i \geq 0$. On the other
hand, $\sum_{i \in Q_0} x_i = 0$, which gives $\bx = 0 \in
\calC_Q$.

Assume now that there is a proper nonempty filter $F$ in $Q$ such
that $\bx_F = 0$. Let $Q'$ and $Q''$ be the full subquivers of $Q$
such that $Q'_0 = F$ and $Q''_0 = Q_0 \setminus F$. Then $\bx =
\bx' + \bx''$ according to the canonical isomorphism $\bbZ^{Q_0}
\simeq \bbZ^{Q'_0} \oplus \bbZ^{Q''_0}$. Observe that $\bx' \in
X_{Q'}$ and $\bx'' \in X_{Q''}$. By the inductive assumption,
$\bx' \in \calC_{Q'}$ and $\bx'' \in \calC_{Q''}$. Consequently,
$\bx \in \calC_{Q'} \oplus \calC_{Q''} \subseteq \calC_Q$.

Hence we may assume that $Q_1$ is nonempty and that $\bx_F > 0$
for any nonempty proper filter $F$ in $Q$. Choose $\alpha \in Q_1$
and let $\by = \bx - \be_\alpha$. Obviously $\by_{Q_0} = 0$. Since
there are no oriented cycles in $Q$, there is a filter $F$ in $Q$
with $t \alpha\in F$ and $s \alpha \not \in F$. For any such
filter $\by_F = \bx_F - 1 \geq 0$, while for the remaining ones $
\by_F = \bx_F \geq 0$. Hence $\by \in X_S$ and $\sum_{F \in \calF}
\by_F < \sum_{F \in \calF} \bx_F$. By our inductive assumption
$\by \in \calC_S$, which gives $\bx = \by + \be_\alpha \in
\calC_S$.
\end{proof}

Now we consider the problem of finding equations defining
$\ol{\calO}_V$. More precisely, we want to describe generators of
the ideal $I_{\calC_Q}$, which is the kernel of the algebra
homomorphism
\[
k [S_\alpha]_{\alpha \in Q_1} \to k [T_i, T_i^{-1}]_{i \in Q_0},
\qquad S_\alpha \mapsto T^{\be_\alpha}.
\]
For $\bw = (w_\alpha)_{\alpha \in Q_1} \in \bbZ^{Q_1}$ we define
$\bw^+ = (w_\alpha^+)_{\alpha \in Q_1}, \bw^- =
(w_\alpha^-)_{\alpha \in Q_1} \in \bbZ^{Q_1}$ by
\[
w^+_\alpha = \max \{ w_\alpha, 0 \} \quad \text {and} \quad
w^-_\alpha = \max \{ -w_\alpha, 0\} \quad \text{for } \alpha \in
Q_1.
\]
Let $\calU : \bbZ^{Q_1} \to \bbZ^{Q_0}$ be the group homomorphism
such that $\calU (\bf_\alpha) = \be_\alpha$ for $\alpha \in Q_1$,
where $(\bf_\alpha)_{\alpha \in Q_1}$ is the standard basis of
$\bbZ^{Q_1}$. Then $I_{\calC_Q}$ is generated by the binomials
$S^{\bw^+} - S^{\bw^-}$ with $\bw \in \Ker (\calU)$, where
$S^{\bw} = \prod_{i \in Q_1} S_\alpha^{w_\alpha}$ for $\bw =
(w_\alpha)_{\alpha \in Q_1} \in \bbN^{Q_1}$ (see \cite{St}*{Lemma
1.1}). Note that $\Ker (\calU)$ consists of the vectors $\bw =
(w_\alpha)_{\alpha \in Q_1} \in \bbZ^{Q_1}$ such that
\begin{equation} \label{balance}
\sum_{s \alpha = i} w_\alpha = \sum_{t \alpha = i} w_\alpha \text{
for all } i \in Q_0.
\end{equation}
In the case of toric varieties occurring in
Theorem~\ref{theotoric} we shall indicate a special finite subsets
of $\Ker (\calU)$ for which the corresponding binomials generate
the ideal $I_{\calC_Q}$.

Let $Q^*$ be the double quiver of $Q$, i.e., the quiver with the
same set of vertices as $Q$ and the set of arrows $Q_1 \cup
Q_1^-$, where $Q_1^- = \{ \alpha^- \mid \alpha \in Q_1 \}$ is the
set of the formal inverses $\alpha^-$ of arrows $\alpha$ in $Q$
with $s \alpha^- = t \alpha$ and $t \alpha^- = s \alpha$. By a
nonoriented path in $Q$ we mean an oriented path in $Q^*$ which
does not contain neither $\alpha \alpha^-$ nor $\alpha^- \alpha$
for $\alpha \in Q_1$ as a subpath. By a nonoriented cycle in $Q$
we mean a nontrivial nonoriented path in $Q$ which starts and
terminates at the same vertex. A nonoriented cycle is called
primitive if it does not contain a proper subpath which is a
nonoriented cycle.

With a primitive nonoriented cycle $\beta_1 \cdots \beta_l$ in $Q$
we may associate a vector $\bu = (u_\alpha)_{\alpha \in Q_1} \in
\bbZ^{Q_1}$ in the following way:
\[
u_\alpha =
\begin{cases}
1, & \alpha = \beta_i \text{ for some } i \in [1, l],
\\ %
-1, & \alpha^- = \beta_i \text{ for some } i \in [1, l],
\\ %
0, & \text{otherwise},
\end{cases} \qquad \alpha \in Q_1.
\]
Note that $\bu \in \Ker (\calU)$. Let $\calZ$ be the set of all
vectors obtained from primitive nonoriented cycles in $Q$ in the
way described above. Observe that $\calZ = - \calZ$, which means
that $-\bu \in \calZ$ for any $\bu \in \calZ$. Thus we can choose
a subset $\calZ'$ of $\calZ$ such that $\calZ = \calZ' \cup
(-\calZ')$ and $\calZ' \cap (-\calZ') = \emptyset$. Note that the
elements of $\calZ'$ correspond bijectively to the equivalence
classes of primitive nonoriented cycles in $Q$ under the relation
which identify a cycle with all its rotations and all rotations of
its inversion (since these notions seem to be self-explained we
will not give precise definitions here). Our next aim is to show
that the binomials corresponding to the elements of $\calZ'$
(hence to the equivalence classes of primitive nonoriented cycles
in $Q$) generate $\Ker (\calU)$. We start with the following
auxiliary observation.

\begin{lemm}
If $\bw \in \Ker (\calU)$ is nonzero, then there exists $\bu \in
\calZ$ such that $\bu^+ \leq \bw^+$ and $\bu^- \leq \bw^-$.
\end{lemm}

\begin{proof}
Let $\bw = (w_\alpha)_{\alpha \in Q_1}$ be a nonzero element of
$\Ker (\calU)$. We construct inductively an infinite nonoriented
path $\omega = \beta_1 \beta_2 \beta_3 \cdots$ in $Q$, such that
for each $j \geq 1$ either $\beta_j = \alpha$ for an arrow $\alpha
\in Q_1$ with $w_\alpha > 0$, or $\beta_j = \alpha^-$ for an arrow
$\alpha \in Q_1$ with $w_\alpha < 0$. We take an arbitrary arrow
$\alpha \in Q_1$ with $w_\alpha \neq 0$ in order to define
$\beta_1$. Assume now that $\beta_n$ is defined. If $\beta_n =
\alpha$ for $\alpha \in Q_1$, then it follows from the
equality~\eqref{balance} for $i = t \alpha_n$ that there is an
arrow $\alpha' \neq \alpha$ such that either $s \alpha' = t
\alpha$ and $w_{\alpha'} > 0$, or $t \alpha' = t \alpha$ and
$w_{\alpha'} < 0$. In the former case we put $\beta_{n + 1} =
\alpha'$, and in the latter $\beta_{n + 1} = \alpha'^-$. If
$\beta_n = \alpha^-$ for $\alpha \in Q_1$, then we consider the
equality~\eqref{balance} for $i = s \alpha$ and we define
$\beta_{n + 1}$ in a similar way as above. Since the quiver $Q$ is
finite, there exists a primitive nonoriented cycle which is a
subpath of $\omega$. The vector corresponding to this cycle
satisfies the claim.
\end{proof}

Now we can prove the announced result.

\begin{prop}\ \label{propideal}
Let $Q$ be a finite quiver without oriented cycles and assume the
above notation. Then the ideal $I_{\calC_Q}$ is generated by the
binomials
\[
S^{\bu^+} - S^{\bu^-}, \quad \bu \in\calZ'.
\]
\end{prop}

\begin{proof}
Since
\[
S^{\bv^+} - S^{\bv^-} =- (S^{\bu^+} - S^{\bu^-})
\]
if $\bv = -\bu$ and $ \bu \in \bbZ^{Q_1}$, it suffices to prove
that if $\bw = (w_\alpha)_{\alpha \in Q_1}$ belongs to
$\Ker(\calU)$, then $S^{\bw^+} - S^{\bw^-}$ belongs to the ideal
generated by the binomials
\[
S^{\bu^+} - S^{\bu^-}, \quad \bu \in \calZ.
\]
We proceed by induction on $|\bw| = \sum_{\alpha \in Q_1}
|w_\alpha| \geq 0$. If $|\bw| = 0$, then $\bw = 0$ and we are
done. Otherwise by the previous lemma, there is a vector $\bu \in
\calZ$ such that $\bu^+ \leq \bw^+$ and $\bu^- \leq \bw^-$. Then
$\bw^+ = \bu^+ + \bv^+$ and $\bw^- = \bu^- + \bv^-$ for  $\bv =
\bw - \bu$. Moreover, $\bv \in \Ker (\calU)$ and $|\bv| = |\bw| -
|\bu| < |\bw|$. Since
\[
S^{\bw^+}-S^{\bw^-} = S^{\bv^+} (S^{\bu^+} - S^{\bu^-}) +
S^{\bu^-} (S^{\bv^+} - S^{\bv^-}),
\]
the claim follows by the inductive assumption.
\end{proof}

The above proposition gives us a finite set of generators of
$I_{\calC_Q}$. As we shall see below, this set usually is not
minimal.

We restrict now our findings to a quiver $Q$ of a special form.
Let $0 \leq p \leq q \leq r \leq s \leq t$. We define a quiver $Q
= Q (p, q, r, s, t)$ in the following way. If $0 < p < q < r < s <
t$, then $Q$ is the quiver
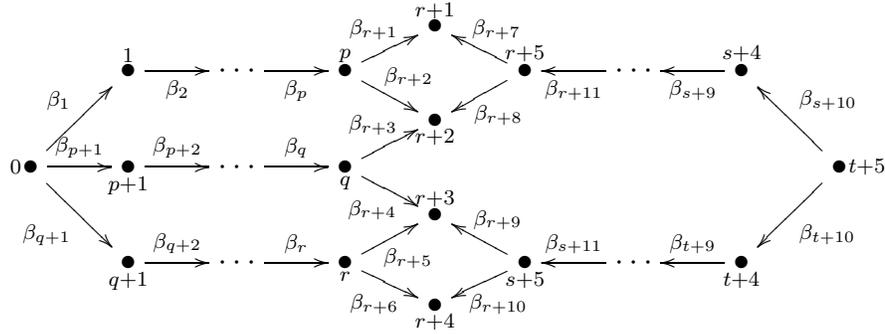
\[
\xymatrix{%
& & & & \ar@{}[d]|(0.5)*+{\textstyle \bullet}="a1"
\ar@{}"a1";[]|<{r + 1} \ar@{}[dd]|(0.75)*+{\textstyle
\bullet}="a2" \ar@{}"a2";[dd]|<{r + 2} & & & &
\\ %
& \bullet \save*+!D{\scriptstyle 1} \restore \ar[r]_-{\beta_2} &
\cdots \ar[r]_-{\beta_p} & \bullet \save*+!D{\scriptstyle p}
\restore \ar"a1"^{\beta_{r + 1}} \ar"a2"^{\beta_{r + 2}} &
\ar@{}[dd]|(0.75)*+{\textstyle \bullet}="a3" \ar@{}"a3";[d]|<{r +
3} & \bullet \save*+!D{\scriptstyle r + 5} \restore
\ar"a1"_{\beta_{r + 7}} \ar"a2"^{\beta_{r + 8}} & \cdots
\ar[l]^-{\beta_{r + 11}} & \bullet \save*+!D{\scriptstyle s + 4}
\restore \ar[l]^-{\beta_{s + 9}}
\\ %
\bullet \save*+!R{\scriptstyle 0} \restore \ar[ru]^{\beta_1}
\ar[r]^{\beta_{p + 1}} \ar[rd]_{\beta_{q + 1}} & \bullet
\save*+!U{\scriptstyle p + 1} \restore \ar[r]^-{\beta_{p + 2}} &
\cdots \ar[r]^-{\beta_q} & \bullet \save*+!U{\scriptstyle q}
\restore \ar"a2"^{\beta_{r + 3}} \ar"a3"_{\beta_{r + 4}} &
\ar@{}[dd]|(0.75)*+{\textstyle \bullet}="a4" \ar@{}"a4";[dd]|<{r +
4} & & & & \bullet \save*+!L{\scriptstyle t + 5} \restore
\ar[lu]_{\beta_{s + 10}} \ar[ld]^{\beta_{t + 10}}
\\ %
& \bullet \save*+!U{\scriptstyle q + 1} \restore \ar[r]^-{\beta_{q
+ 2}} & \cdots \ar[r]^-{\beta_r} & \bullet \save*+!U{\scriptstyle
r} \restore \ar"a3"_{\beta_{r + 5}} \ar"a4"_{\beta_{r + 6}} & &
\bullet \save*+!U{\scriptstyle s + 5} \restore \ar"a3"_{\beta_{r +
9}} \ar"a4"^{\beta_{r + 10}} & \cdots \ar[l]_-{\beta_{s + 11}} &
\bullet \save*+!U{\scriptstyle t + 4} \restore \ar[l]_-{\beta_{t +
9}}
\\ %
& & & & & & & &}
\]
If $0 = p$ ($p = q$, $q = r$, $r = s$ or $s = t$, respectively)
then we cancel appropriate arrows and identify vertices $0$ and
$p$ ($0$ and $q$, $0$ and $r$, $r + 5$ and $t + 5$, or $s + 5$ and
$t + 5$, respectively). Thus in the most extremal case $0 = p = q
= r = s = t$ we get the quiver
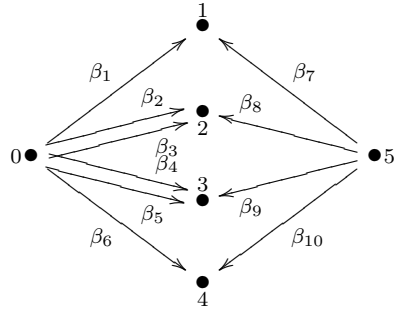
\[
\xymatrix{%
& & \ar@{}[d]|(0.5)*+{\textstyle \bullet}="a1" \ar@{}"a1";[]|<1 &
&
\\ %
& & \ar@{}[d]|(0.5)*+{\textstyle \bullet}="a2" \ar@{}"a2";[d]|<2
\ar@{}[dd]|(0.75)*+{\textstyle \bullet}="a3" \ar@{}"a3";[d]|<3
\ar@{}[ddd]|(0.83)*+{\textstyle \bullet}="a4" \ar@{}"a4";[ddd]|<4
& &
\\ %
\bullet \save*+!R{\scriptstyle 0} \restore \ar"a1"^{\beta_1}
\ar@<0.5ex>"a2"^(0.75){\beta_2} \ar@<-0.5ex>"a2"_(0.75){\beta_3}
\ar@<0.5ex>"a3"^(0.75){\beta_4} \ar@<-0.5ex>"a3"_(0.75){\beta_5}
\ar"a4"_{\beta_6} & & & & \bullet \save*+!L{\scriptstyle 5}
\restore \ar"a1"_{\beta_7} \ar"a2"_(0.75){\beta_8}
\ar"a3"^(0.75){\beta_9} \ar"a4"^{\beta_{10}}
\\ %
& & & &
\\ %
& & & &}
\]
with $6$ vertices and $10$ arrows.

Recall that $\bf_{\beta_1}$, \ldots, $\bf_{\beta_{t + 10}}$ is the
standard basis of $\bbZ^{Q_1}$. Let $\bu_i = \bf_{\beta_{r + i}}$
for $i \in [1, 10]$ and
\begin{gather*}
\bu_{11} = \bf_{[1, p]}, \; \bu_{12} =\bf_{[p + 1, q]},\; \bu_{13}
= \bf_{[q + 1, r]},
\\ %
\bu_{14} = \bf_{[r + 11, s + 10]}, \; \bu_{15} = \bf_{[s + 11, t +
10]},
\end{gather*}
where $\bf_{[i, j]} = \sum_{l \in [i, j]} \bf_{\beta_l}$ for $i, j
\in [1, t + 10]$. Observe that it may happen that $\bu_i = 0$ for
some $i \in [11, 15]$. With the above notation $\calZ'$ consists,
up to sign, of the following vectors:
\begin{align*}
\bv_1 & = \bu_2 + \bu_{11} - \bu_3 - \bu_{12},
\\ %
\bv_2 & = \bu_4 + \bu_{12} - \bu_5 - \bu_{13},
\\ %
\bv_3 & = \bu_1 + \bu_8 - \bu_2 - \bu_7,
 \\ %
\bv_4 & = \bu_5 + \bu_{10} - \bu_6 - \bu_9,
\\ %
\bv_5 & = \bu_3 + \bu_9 + \bu_{15} - \bu_4 - \bu_8 - \bu_{14},
\\ %
\bv_6 & = \bu_1 + \bu_9 + \bu_{11} + \bu_{15} - \bu_4 - \bu_7 -
\bu_{12} - \bu_{14},
\\ %
\bv_7 & = \bu_3 + \bu_{10} + \bu_{12} + \bu_{15} - \bu_6 - \bu_8 -
\bu_{13} - \bu_{14},
\\ %
\bv_8 & = \bu_1 + \bu_{10} + \bu_{11} + \bu_{15} - \bu_6 - \bu_7 -
\bu_{13} - \bu_{14},
\\ %
\bv_9 & = \bu_1 + \bu_8 + \bu_{11} - \bu_3 - \bu_7 - \bu_{12},
\\ %
\bv_{10} & = \bu_4 + \bu_{10} + \bu_{12} - \bu_6 - \bu_9 - \bu_{13},
\\ %
\bv_{11} & = \bu_2 + \bu_4 + \bu_{11} - \bu_3 - \bu_5 - \bu_{13},
\\ %
\bv_{12} & = \bu_1 + \bu_3 + \bu_9 + \bu_{15} - \bu_2 - \bu_4 -
\bu_7 - \bu_{14},
\\ %
\bv_{13} & = \bu_3 + \bu_5 + \bu_{10} + \bu_{15} - \bu_4 - \bu_6 -
\bu_8 - \bu_{14},
\\ %
\bv_{14} & = \bu_2 + \bu_9 + \bu_{11} + \bu_{15} - \bu_4 - \bu_8 -
\bu_{12} - \bu_{14},
\\ %
\bv_{15} & = \bu_3 + \bu_9 + \bu_{12} + \bu_{15} - \bu_5 - \bu_8 -
\bu_{13} - \bu_{14},
\\ %
\bv_{16} & = \bu_1 + \bu_4 + \bu_8 + \bu_{11} - \bu_3 - \bu_5 -
\bu_7 - \bu_{13},
\\ %
\bv_{17} & = \bu_2 + \bu_4 + \bu_{10} + \bu_{11} - \bu_3 - \bu_6 -
\bu_9 - \bu_{13},
\\ %
\bv_{18} & = \bu_1 + \bu_3 + \bu_9 + \bu_{12} + \bu_{15} - \bu_2 -
\bu_5 - \bu_7 - \bu_{13} - \bu_{14},
\\ %
\bv_{19} & = \bu_2 + \bu_5 + \bu_{10} + \bu_{11} + \bu_{15} -
\bu_4 - \bu_6 - \bu_8 - \bu_{12} - \bu_{14},
\\ %
\bv_{20} & = \bu_1 + \bu_3 + \bu_5 + \bu_{10} + \bu_{15} - \bu_2 -
\bu_4 - \bu_6 - \bu_7 - \bu_{14},
\\ %
\bv_{21} & = \bu_2 + \bu_9 + \bu_{11} + \bu_{15} - \bu_5 - \bu_8 -
\bu_{13} - \bu_{14},
\\ %
\bv_{22} & = \bu_1 + \bu_4 + \bu_8 + \bu_{10} + \bu_{11} - \bu_3 -
\bu_6 - \bu_7 - \bu_9 - \bu_{13},
\\ %
\bv_{23} & = \bu_1 + \bu_9 + \bu_{11} + \bu_{15} - \bu_5 - \bu_7 -
\bu_{13} - \bu_{14},
\\ %
\bv_{24} & = \bu_2 + \bu_{10} + \bu_{11} + \bu_{15} - \bu_6 -
\bu_8 - \bu_{13} - \bu_{14},
\\ %
\bv_{25} & = \bu_1 + \bu_5 + \bu_{10} + \bu_{11} + \bu_{15} -
\bu_4 - \bu_6 - \bu_7 - \bu_{12} - \bu_{14},
\\ %
\bv_{26} & = \bu_1 + \bu_3 + \bu_{10} + \bu_{12} + \bu_{15} -
\bu_2 - \bu_6 - \bu_7 - \bu_{13} - \bu_{14}.
\end{align*}
Indeed, recall that the elements of $\calZ'$ correspond to the
equivalence classes of the primitive nonoriented cycles in $Q$.
Note that each such equivalence class is determined by a nonempty
subset of the set consisting of the five inner polygons visible on
the picture of the quiver $Q$. There are $2^5 - 1 = 31$ such
nonempty subsets, $26$ of them leads to our vectors $\bv_i$, $i
\in [1, 26]$, and none of the remaining five subsets corresponds
to the equivalence class of a primitive nonoriented cycle in $Q$
(they may be seen as corresponding to equivalence classes of two
disjoint primitive cycles).

\begin{lemm} \label{lemmideal}
Let $Q = Q (p, q, r, s, t)$ for $0 \leq p \leq q \leq r \leq s
\leq t$. Then the ideal $I_{\calC_Q}$ is generated by the
binomials
\[
S^{\bv_i^+} - S^{\bv_i^-}, \quad i \in [1,8].
\]
\end{lemm}

\begin{proof}
By Proposition~\ref{propideal}, it suffices to show that the above
binomials generate the remaining binomials
$$
S^{\bv_i^+} - S^{\bv_i^-}, \quad i \in [9,26].
$$
This is a quite easy, but tedious verification. Hence we prove the
claim only for $i = 9$ and $i = 21$, leaving the other cases to
the reader:
\begin{align*}
S^{\bv_9^+} -  S^{\bv_9^-} & = S^{\bu_1} S^{\bu_8} S^{\bu_{11}} -
S^{\bu_3} S^{\bu_7} S^{\bu_{12}}
\\ %
& = S^{\bu_{11}} (S^{\bu_1} S^{\bu_8} - S^{\bu_2} S^{\bu_7}) +
S^{\bu_7}  (S^{\bu_2} S^{\bu_{11}} - S^{\bu_3} S^{\bu_{12}})
\\ %
& = S^{\bu_{11}} (S^{\bv_3^+} - S^{\bv_3^-}) + S^{\bu_7}
(S^{\bv_1^+} - S^{\bv_1^-}),
\\ %
S^{\bv_{21}^+} - S^{\bv_{21}^-} & = S^{\bu_2} S^{\bu_9}
S^{\bu_{11}} S^{\bu_{15}} - S^{\bu_5} S^{\bu_8} S^{\bu_{13}}
S^{\bu_{14}}
\\ %
& = S^{\bu_9} S^{\bu_{15}} (S^{\bu_2} S^{\bu_{11}} - S^{\bu_3}
S^{\bu_{12}})
\\ %
& \quad + S^{\bu_{12}} (S^{\bu_3} S^{\bu_9} S^{\bu_{15}}-
S^{\bu_4} S^{\bu_8} S^{\bu_{14}})
\\ %
& \quad + S^{\bu_8} S^{\bu_{14}} (S^{\bu_4} S^{\bu_{12}} -
S^{\bu_5} S^{\bu_{13}})
\\ %
& = S^{\bu_9} S^{\bu_{15}} (S^{\bv_1^+} - S^{\bv_1^-}) +
S^{\bu_{12}} (S^{\bv_5^+} - S^{\bv_5^-})
\\ %
& \quad + S^{\bu_8} S^{\bu_{14}} (S^{\bv_2^+} - S^{\bv_2^-}).
\qedhere
\end{align*}
\end{proof}

\section{Deformations to toric varieties} \label{sect3}

Let $\Delta$, $\bd$ and $P$ be as in Theorem~\ref{theospec}. As
usual $\be_1$, \ldots, $\be_{t + 5}$ denote the standard basis of
$\bbZ^{t + 5}$. For $i, j \in [1, t + 5]$, $\be_{[i, j]} = \sum_{l
\in [i, j]} \be_l$. If $\bx = (x_i)_{i \in [1, t + 5]} \in k^{t +
5}$ and $\bw = (w_i)_{i \in [1, t + 5]} \in \bbN^{t + 5}$, then
$\bx^{\bw} = \prod_{i \in [1, t+ 5]} x_i^{w_i}$.

Our aim in this section is to prove Theorem~\ref{theospec}. As the
first step we describe the coordinate ring of $\ol{\calO}_P$. Note
that $\dim \ol{\calO}_P = t + 5$. Indeed, $\dim \ol{\calO}_P =
\dim \GL (\bd) - \dim \Stab_{\GL (\bd)} (P)$, where $\Stab_{\GL
(\bd)}$ denotes the subgroup of all $g \in \GL (\bd)$ such that $g
\cdot P = P$. Easy calculations show $\dim \GL (\bd) = t + 6$ and
$\Stab_{\GL (\bd)} (P) \simeq k^*$, thus the formula follows.

Let $\Phi : k^{t + 5} \to \rep_\Delta (\bd)$ be given by
\begin{align*}
\Phi (\bx)_{\alpha_i} & = [x_i], \; i \in [1, r] \cup [r + 6, t +
5],
\\ %
\Phi (\bx)_{\alpha_{r + 1}} & = \bx^{\be_{[p + 1, r]}} [
\begin{matrix}
x_{r + 1} & x_{r + 3}
\end{matrix}
],
\\ %
\Phi (\bx)_{\alpha_{r + 2}} & = \bx^{\be_{[1, p]}} \bx^{\be_{[q +
1, r]}} [
\begin{matrix}
- x_{r + 1} - x_{r + 4} & - x_{r + 2} - x_{r + 3}
\end{matrix}
],
\\ %
\Phi (\bx)_{\alpha_{r + 3}} & = \bx^{\be_{[1, q]}} [
\begin{matrix}
x_{r + 4} & x_{r + 2}
\end{matrix}
],
\\ %
\Phi (\bx)_{\alpha_{r + 4}} & = [
\begin{matrix}
- x_{r + 3} & x_{r + 1}
\end{matrix}
]^{\tr} x_{r + 5} \bx^{\be_{[s + 6, t + 5]}},
\\ %
\Phi (\bx)_{\alpha_{r + 5}} & = [
\begin{matrix}
x_{r + 2} & - x_{r + 4}
\end{matrix}
]^{\tr} x_{r + 5} \bx^{\be_{[r + 6, s + 5]}},
\end{align*}
for $\bx = (x_i)_{i \in [1, t + 5]} \in k^{t + 5}$. The next
observation is the following.

\begin{lemm} \label{lemmcoord}
$\ol{\Phi (k^{t + 5})} = \ol{\calO}_P$.
\end{lemm}

\begin{proof}
Let
\begin{multline*}
U = \{ \bx = (x_i)_{i \in [1, t + 5]} \in k^{t + 5} \mid x_i \neq
0, \, i \in [1, r] \cup [r + 5, t + 5],
\\ %
x_{r + 1} x_{r + 2} \neq x_{r + 3} x_{r + 4} \}.
\end{multline*}
Then $U$ is an open subset of $k^{t + 5}$ and $\Phi |_U$ is
injective, thus we get $\dim \ol{\Phi (k^{t + 5})} = t + 5 = \dim
\ol{\calO}_P$. Since $\ol{\calO}_P$ is irreducible, it is enough
to show that $\Phi (U) \subset \calO_P$.
Let $\bx = (x_i)_{i \in [1, t + 5]} \in U$ and $X = \left[
\begin{smallmatrix}
x_{r + 1} & x_{r + 3} \\ x_{r + 4} & x_{r + 2}
\end{smallmatrix}
\right]$. Then $g = (g_i)_{i \in [1, t + 2]}$ given by
\begin{align*}
g_i & = \bx^{\be_{[1, i]}}, \; i \in [0, p],
\\ %
g_i & = \bx^{\be_{[p + 1, i]}}, \; i \in [p + 1, q],
\\ %
g_i & = \bx^{\be_{[q + 1, i]}}, \; i \in [q + 1, r],
\\ %
g_{r + 1} & = \bx^{\be_{[1, r]}} X,
\\ %
g_i & = \bx^{\be_{[1, r]}} \det X x_{r + 5} \bx^{\be_{[r + 6, i +
3]}} \bx^{\be_{[s + 6, t + 5]}}, \; i \in [r + 2, s + 1],
\\ %
g_i & = \bx^{\be_{[1, r]}} \det X x_{r + 5} \bx^{\be_{[r + 6, s +
5]}} \bx^{\be_{[s + 6, i + 3]}}, \; i \in [s + 2, t + 2],
\end{align*}
belongs to $\GL (\bd)$ and $g \cdot \Phi (\bx) = P$.
\end{proof}

An obvious reformulation of the above lemma says that $k
[\ol{\calO}_P] = k [a_1, \ldots, a_{t + 10}]$, where $a_1, \ldots,
a_{t + 10}$ are polynomials in $k [T_1, \ldots, T_{t + 5}]$
defined by
\begin{align*}
a_i & = T_i, \; i \in [1, r],
\\ %
a_{r + 1} & = T^{\be_{[p + 1, r]}} T_{r + 1},
\\ %
a_{r + 2} & = T^{\be_{[p + 1, r]}} T_{r + 3},
\\ %
a_{r + 3} & = T^{\be_{[1, p]}} T^{\be_{[q + 1, r]}} T_{r + 2} +
T^{\be_{[1, p]}} T^{\be_{[q + 1, r]}} T_{r + 3},
\\ %
a_{r + 4} & = T^{\be_{[1, p]}} T^{\be_{[q + 1, r]}} T_{r + 1} +
T^{\be_{[1, p]}} T^{\be_{[q + 1, r]}} T_{r + 4},
\\ %
a_{r + 5} & = T^{\be_{[1, q]}} T_{r + 4},
\\ %
a_{r + 6} & = T^{\be_{[1, q]}} T_{r + 2},
\\ %
a_{r + 7} & = T_{r + 1} T_{r + 5} T^{\be_{[s + 6, t + 5]}},
\\ %
a_{r + 8} & = T_{r + 3} T_{r + 5} T^{\be_{[s + 6, t + 5]}},
\\ %
a_{r + 9} & = T_{r + 4} T_{r + 5} T^{\be_{[r + 6, s + 5]}},
\\ %
a_{r + 10} & = T_{r + 2} T_{r + 5} T^{\be_{[r + 6, s + 5]}},
\\ %
a_i & = T_{i - 5}, \; i \in [r + 11, t + 10].
\end{align*}
As before, $T^{\bw} = \prod_{i \in [1, t + 10]} T_i^{w_i}$ for
$\bw = (w_i)_{i \in [1, t + 10]} \in \bbN^{t + 10}$.

We order the elements of $\bbN^{t + 5}$ by the reversed
lexicographic order, i.e., we say that $\bu = (u_i)_{i \in [1, t +
5]}$ is smaller than $\bv = (v_i)_{i \in [1, t + 5]}$ if there
exists $i \in [1, t + 5]$ such that $u_i < v_i$ and $u_j = v_j$
for all $j \in [i + 1, t + 5]$. The induced order of the monomials
in $k [T_1, \ldots, T_{t + 5}]$ is a term order in the sense
of~\cite{RoSw}*{1.3}.

For $a = \sum_{\bv \in \bbN^{t + 5}} \lambda_{\bv} T^{\bv} \in k
[T_1, \ldots T_{t + 5}]$, $a \neq 0$, we define the initial
monomial $\ini (a)$ as $T^{\bu}$, where $\bu = \max \{ \bv \in
\bbN^{t + 5} \mid \lambda_{\bv} \neq 0 \}$. If $A$ is a subalgebra
of $k [T_1, \ldots, T_{t + 5}]$, then by the initial algebra $\ini
(A)$ of $A$ we mean the subalgebra of $A$ generated by $\{ \ini
(a) \mid a \in A \}$. According
to~\cite{CoHeVa}*{Corollary~2.3(b)} in order to prove
Theorem~\ref{theospec} it is enough to show that $\ini (k [a_1,
\ldots, a_{t + 10}])$ is finitely generated and normal. Using
Theorem~\ref{theotoric} it will follow if we show isomorphisms
$\ini (k [a_1, \ldots, a_{t + 10}]) \simeq k [\ini (a_1), \ldots,
\ini (a_{t + 10})] \simeq k [\ol{\calO}_V]$, where $V$ is the
point of $\rep_Q ((1)_{i \in [1, t + 5]})$ with all matrices equal
to $[1]$. Here $Q = Q (p, q, r, s, t)$ is the quiver defined in
Section~\ref{sect2}.

We first show the latter isomorphism, or in other words, we
describe $k [\ol{\calO}_V]$. The method is analogous to the one
applied above in order to describe $k [\ol{\calO}_P]$. Let $\Psi :
k^{t + 5} \to \rep_Q ((1)_{i \in [1, t + 5]})$ be defined by
\begin{align*}
\Phi (\bx)_{\beta_i} & = x_i, \; i \in [1, r],
\\ %
\Phi (\bx)_{\beta_{r + 1}} & = \bx^{\be_{[p + 1, r]}} x_{r + 1},
\\ %
\Phi (\bx)_{\beta_{r + 2}} & = \bx^{\be_{[p + 1, r]}} x_{r + 3},
\\ %
\Phi (\bx)_{\beta_i} & = \bx^{\be_{[1, p]}} \bx^{\be_{[q + 1, r]}}
x_i, \; i \in [r + 3, r + 4],
\\ %
\Phi (\bx)_{\beta_{r + 5}} & = \bx^{\be_{[1, q]}} x_{r + 4},
\\ %
\Phi (\bx)_{\beta_{r + 6}} & = \bx^{\be_{[1, q]}} x_{r + 2},
\\ %
\Phi (\bx)_{\beta_{r + 7}} & = x_{r + 1} x_{r + 5} \bx^{\be_{[s +
6, t + 5]}},
\\ %
\Phi (\bx)_{\beta_{r + 8}} & = x_{r + 3} x_{r + 5} \bx^{\be_{[s +
6, t + 5]}},
\\ %
\Phi (\bx)_{\beta_{r + 9}} & = x_{r + 4} x_{r + 5} \bx^{\be_{[r +
6, s + 5]}},
\\ %
\Phi (\bx)_{\beta_{r + 10}} & = x_{r + 2} x_{r + 5} \bx^{\be_{[r +
6, s + 5]}},
\\ %
\Phi (\bx)_{\beta_i} & = x_{i - 5}, \; i \in [r + 11, t + 10],
\end{align*}
for $\bx = (x_i)_{i \in [1, t + 5]} \in k^{t + 5}$. With arguments
similar to those used in the proof of Lemma~\ref{lemmcoord}, one
shows that $\ol{\Phi (k^{t + 5})} = \ol{\calO}_V$, hence $k
[\ol{\calO}_V]$ may be identified with the subalgebra of $k [T_1,
\ldots, T_{t + 5}]$ generated by polynomials $b_1$, \ldots, $b_{t
+ 10}$, where
\begin{align*}
b_i & = T_i, \; i \in [1, r],
\\ %
b_{r + 1} & = T^{\be_{[p + 1, r]}} T_{r + 1},
\\ %
b_{r + 2} & = T^{\be_{[p + 1, r]}} T_{r + 3},
\\ %
b_i & = T^{\be_{[1, p]}} T^{\be_{[q + 1, r]}} T_i, \; i \in [r +
3, r + 4],
\\ %
b_{r + 5} & = T^{\be_{[1, q]}} T_{r + 4},
\\ %
b_{r + 6} & = T^{\be_{[1, q]}} T_{r + 2},
\\ %
b_{r + 7} & = T_{r + 1} T_{r + 5} T^{\be_{[s + 6, t + 5]}},
\\ %
b_{r + 8} & = T_{r + 3} T_{r + 5} T^{\be_{[s + 6, t + 5]}},
\\ %
b_{r + 9} & = T_{r + 4} T_{r + 5} T^{\be_{[r + 6, s + 5]}},
\\ %
b_{r + 10} & = T_{r + 2} T_{r + 5} T^{\be_{[r + 6, s + 5]}},
\\ %
b_i & = T_{i - 5}, \; i \in [r + 11, t + 10].
\end{align*}
It is an obvious observation that $b_i = \ini (a_i)$ for all $i
\in [1, t + 10]$, which shows that $k [\ini (a_1), \ldots, \ini
(a_{t + 10})] \simeq k [\ol{\calO}_V]$.

Observe that the kernel $I$ of the algebra homomorphism
\[
k [S_{\beta_1}, \ldots, S_{\beta_{t + 10}}] \to k [T_1, \ldots,
T_{t + 5}], \qquad S_{\beta_i} \mapsto b_i,
\]
equals the ideal $I_{C_Q}$ defined in Section~\ref{sect2},
as both of them are the ideals of $\ol{\calO}_V$ in
$\rep_Q ((1)_{i \in [1, t + 5]})$. By Lemma~\ref{lemmideal},
$I$ is generated by the binomials $\xi_i = S^{\bv_i^+} -
S^{\bv_i^-}$, $i \in [1, 8]$, where $\bv_1$, \ldots, $\bv_8$
are as in Section~\ref{sect2}.

As the final step we show that $\ini (k [a_1, \ldots, a_{t + 10}])
\simeq k [b_1, \ldots, b_{t + 10}]$ (if this condition holds, then
one says that $a = ( a_1, \ldots, a_{t + 10} )$ is a Sagbi basis
of the algebra $k [a_1, \ldots, a_{t + 10}]$). According
to~\cite{CoHeVa}*{Proposition~1.1} it is enough to show that there
exist $\lambda_{i, \bu} \in k$, $i \in [1, 8]$, $\bu \in I_i = \{
\bv \in \bbN^{t + 10} \mid \ini (a^{\bv}) \leq \ini (\xi_i (a))
\}$, such that
\[
\xi_i (a) = \sum_{\bu \in I_i} \lambda_{i, \bu} a^{\bu}.
\]
Here, $a^{\bu} = a_1^{u_{\beta_1}} \cdots a_{t + 10}^{u_{\beta_{t
+ 10}}}$ for $\bu = (u_{\beta_i})_{i \in [1, t + 10]} \in
\bbN^{Q_1}$ and, for $\xi \in k [S_{\beta_1}, \ldots, S_{\beta_{t
+ 10}}]$, $\xi (a)$ denotes the image of $\xi$ via the map
\[
k [S_{\beta_1}, \ldots, S_{\beta_{t + 10}}] \to k [T_1, \ldots,
T_{t + 5}], \qquad S_{\beta_i} \mapsto a_i.
\]
But
\begin{align*}
\xi_i (a) & = 0, \, i \in \{ 3, 4, 8 \},
\\ %
\xi_1 (a) & = - T^{\be_{[1, r]}} T_{r + 2} = - a^{\be_{[q + 1, r]}}
a_{r + 6},
\\ %
\xi_2 (a) & = T^{\be_{[1, r]}} T_{r + 1} = a^{\be_{[1, p]}} a_{r +
1},
\\ %
\xi_6 (a) & = - T^{\be_{[1, r]}} T_{r + 1} T_{r + 1} T_{r + 5}
T^{\be_{[r + 6, t + 5]}}
\\ %
& = - a^{\be_{[1, p]}} a_{r + 1} a_{r + 7} a^{\be_{[r + 11, s + 10]}},
\\ %
\xi_7 (a) & = T^{\be_{[1, r]}} T_{r + 2} T_{r + 2} T_{r + 5}
T^{\be_{[r + 6, t + 5]}}
\\ %
& = a^{\be_{[q + 1, r]}} a_{r + 6} a_{r + 10} a^{\be_{[s + 11,
t + 10]}},
\\ %
\xi_5 (a) & = T^{\be_{[1, p]}} T^{\be_{[q + 1, r]}} T_{r + 2} T_{r
+ 4} T_{r + 5} T^{\be_{[r + 6, t + 5]}}
\\ %
& \qquad - T^{\be_{[1, p]}} T^{\be_{[q + 1, r]}} T_{r + 1} T_{r +
3} T_{r + 5} T^{\be_{[r + 6, t + 5]}}
\\ %
& = a_{r + 4} a_{r + 10} a^{\be_{[s + 11, t + 10]}} - a_{r + 3}
a_{r + 7} a^{\be_{[r + 11, s + 10]}},
\end{align*}
and the initial monomial
\[
\ini (a_{r + 3} a_{r + 7} a^{\be_{[r + 11, s + 10]}}) =
T^{\be_{[1, p]}} T^{\be_{[q + 1, r]}} T_{r + 1} T_{r + 3} T_{r +
5} T^{\be_{[r + 6, t + 5]}}
\]
is smaller than
\begin{align*}
\ini (a_{r + 4} a_{r + 10} a^{\be_{[s + 11, t + 10]}} ) & = \ini
(\xi_5 (a)) \\ %
& = T^{\be_{[1, p]}} T^{\be_{[q + 1, r]}} T_{r + 2} T_{r + 4} T_{r
+ 5} T^{\be_{[r + 6, t + 5]}},
\end{align*}
which finishes the proof.

\end{document}